\documentclass[12pt]{amsart}
\usepackage{amsmath,amssymb}

\newcommand{\Z}{{\mathbb{Z}}}
\newcommand{\Q}{{\mathbb{Q}}}
\newcommand{\F}{{\mathbb{F}}}
\newcommand{\C}{{\mathbb{C}}}
\newcommand{\ba}{{\mathbf{a}}}
\newcommand{\bc}{{\mathbf{c}}}
\newcommand{\cD}{{\mathcal{D}}}
\newcommand{\cH}{{\mathcal{H}}}
\newcommand{\cL}{{\mathcal{L}}}
\newcommand{\cR}{{\mathcal{R}}}

\newcommand{\fs}{{\mathfrak{s}}}
\newcommand{\ft}{{\mathfrak{t}}}
\newcommand{\fu}{{\mathfrak{u}}}
\newcommand{\fv}{{\mathfrak{v}}}
\newcommand{\Irr}{{\operatorname{Irr}}}

\renewcommand{\leq}{\leqslant}
\renewcommand{\geq}{\geqslant}
\renewcommand{\atop}[2]{\genfrac{}{}{0pt}{}{#1}{#2}}
\newtheorem{theorem}{Theorem}[section]

\newtheorem{corollary}[theorem]{Corollary}
\newtheorem{proposition}[theorem]{Proposition}

\theoremstyle{definition}
\newtheorem{example}[theorem]{Example}

\theoremstyle{remark}
\newtheorem{remark}[theorem]{Remark}

\begin{document}
\title{Hecke algebras of finite type are cellular}
\author{Meinolf Geck}
%\institute{Department of Mathematical Sciences, King's College,
%Aberdeen University, Aberdeen AB24 3UE, Scotland, U.K.,
%\email{m.geck@maths.abdn.ac.uk}}
\address{Department of Mathematical Sciences, King's College,
Aberdeen University, Aberdeen AB24 3UE, Scotland, U.K.}
\email{m.geck@maths.abdn.ac.uk}
\date{March, 2007}
%\maketitle
%\subjclass[2000]{Primary 20C08; Secondary 20C33}

\begin{abstract}
Let $\cH$ be the one-parameter Hecke algebra associated to a finite Weyl 
group $W$, defined over a ground ring in which ``bad'' primes for $W$ are 
invertible. Using deep properties of the Kazhdan--Lusztig basis of $\cH$ 
and Lusztig's $\ba$-function, we show that $\cH$ has a natural cellular 
structure in the sense of Graham and Lehrer. Thus, we obtain a general 
theory of ``Specht modules'' for Hecke algebras of finite type. Previously, 
a general cellular structure was only known to exist in types $A_n$ and $B_n$. 
\end{abstract}

\maketitle
%\pagestyle{myheadings}
%\markboth{Geck}{Hecke algebras are cellular}

%%%%%%%%%%%%%%%%%%%%%%%%%%%%%%%%%%%%%%%%%%%%%%%%%%%%%%%%%%%%%%%%%%%%%%%%%%%%%
\section{Introduction} \label{sec:intro}

The concept of ``cellular algebras'' was introduced by Graham and Lehrer
\cite{GrLe}. It provides a systematic framework for studying the 
representation theory of non-semisimple algebras which are deformations of 
semisimple ones. The original definition was modeled on properties of
the Kazhdan--Lusztig basis \cite{KaLu} in Hecke algebras of type $A_n$.
There is now a significant literature on the subject, and many classes 
of algebras have been shown to admit a ``cellular'' structure, including 
Ariki-Koiki algebras, $q$-Schur algebras, Temperley--Lieb algebras, 
and a variety of other algebras with geometric connections; see, e.g., 
\cite{GrLe}, \cite{GrLe2}, \cite{DJMa}, \cite{CKoe} and the references 
there. However, the question of whether all Hecke algebras of finite type 
(i.e.,  the originally motivating examples) are cellular remained open. 
A positive answer to this question would provide a general theory of 
``Specht modules'' which so far has only been established in types $A_n$, 
$B_n$ (see \cite{DJ0}, \cite{Mu2}, \cite{DJM3}, \cite{GrLe}) and, with 
some restrictions on the ground ring, in type $D_n$ with $n$ odd (see 
\cite{pall1}, \cite{kor2}). In these cases, the constructions heavily rely
on the underlying combinatorics of Young tableaux.

The purpose of this paper is to solve this problem in general. We prove:

\begin{theorem} \label{thmb} Let $\cH$ be the one-parameter Hecke algebra 
associated with a finite Weyl group $W$, defined over an integral domain
in which all ``bad'' primes for $W$ are invertible. Then $\cH$ admits a 
natural ``cellular'' structure, where the elements in the ``cellular
basis'' are integral linear combinations of Kazhdan--Lusztig basis
elements $C_w$ with constant value $\ba(w)$ (Lusztig's $\ba$-function
\cite{Lu1}).
\end{theorem}

If we apply this in type $A_n$, the linear combinations will have only 
one non-zero term and we recover the Kazhdan--Lusztig basis in this case.
In type $B_n$, we obtain a cellular structure which is different from
the one by Dipper--James--Murphy \cite{DJM3} or Graham--Lehrer \cite{GrLe}.
Our construction works in the general setting of multi-parameter Hecke 
algebras, assuming that Lusztig's conjectures on Hecke algebras with 
unequal parameters in \cite[Chapter~14]{Lusztig03} hold. We can also 
formulate a similar result for Hecke algebras associated to 
non-crystallographic finite Coxeter groups $W$; some additional care is 
needed since $\Q$ is not a splitting field for such a group $W$.
%For appropriate
%choices of the ground ring, we can even extend this to the case where
%$W$ is any finite Coxeter group, and not necessarily a Weyl group.

Thus, if ``bad'' primes are invertible in the ground ring, then all Hecke 
algebras of finite type  are ``cellular''. The general theory of cellular 
algebras then produces ``cell representations'' and a natural 
parametrisation of the irreducible representations for non-semisimple 
versions of $\cH$, in terms of properties of certain bilinear forms on the 
cell representations; see \cite[\S 3]{GrLe}. We show that this 
parametrization is precisely given by the ``canonical basic sets'' 
introduced by Rouquier and the author \cite{mykl}, \cite{GeRo2}, 
\cite{mylaus}.  

This paper is organised as follows. In Section~2 we recall the fundamental
facts from Kazhdan--Lusztig theory, including Lusztig's construction of 
the ``asymptotic ring'' $J$. We also establish a basic result about the
representations of $J$ in Proposition~\ref{lem1}. In Section~3, we prove 
Theorem~\ref{thmb} and briefly discuss non-crystallographic Coxeter groups 
in Remark~\ref{noncryst}. The cellular basis that we construct turns out to 
be unique up to integral equivalence of the representations of $J$. Finally, 
in Section~4, we consider examples and discuss the applications to
non-semisimple specialisations. For types $A_n$ and $B_n$, we explain the 
relation of our cellular bases with the previously known ones.

%%%%%%%%%%%%%%%%%%%%%%%%%%%%%%%%%%%%%%%%%%%%%%%%%%%%%%%%%%%%%%%%%%%%%%%%%%%%%
\section{Kazhdan--Lusztig theory and the asymptotic ring} \label{sec2}
References for this section are the books \cite{GePf} and \cite{Lusztig03}.
Let $W$ be a finite Weyl group, with generating set $S$. Let $L \colon 
W \rightarrow \Z$ be a weight function. Thus, we have $L(ww')=L(w)+L(w')$ 
whenever $l(ww')=l(w)+l(w')$ for $w,w'\in W$; here $l(w)$ denotes the length
of $w \in W$. Note that $L(s)=L(t)$ whenever $s,t \in S$ are conjugate in $W$.

We will assume throughout that $L(s)\geq 0$ for all $s \in S$.

Let $R \subseteq \C$ be a subring and $A=R[v,v^{-1}]$ be the ring of 
Laurent polynomials in an indeterminate~$v$. Let $\cH=\cH_A(W,L)$ be
the corresponding generic {\em Iwahori--Hecke algebra}. This is an 
associative algebra, free over $A$ with a basis $\{T_w \mid w \in W\}$ 
and multiplication given as follows.
\[ T_sT_w=\left\{\begin{array}{cl} T_{sw} & \quad \mbox{if $l(sw)>l(w)$},\\
T_{sw}+(v^{L(s)}-v^{-L(s)})T_w & \quad \mbox{if $l(sw)<l(w)$},\end{array}
\right.\]
where $s\in S$ and $w\in W$. 

\subsection{Irreducible representations of $W$ and $\cH$} \label{sub21} 

It is known that $\Q$ is a splitting field for $W$; see, for example,
\cite[6.3.8]{GePf}. We will write
\[ \Irr(W)=\{E^\lambda \mid \lambda \in \Lambda\}, \qquad  d_\lambda= 
\dim E^\lambda,\]
for the set of irreducible representations of $W$ (up to equivalence),
where $\Lambda$ is some finite indexing set. 

Now let $K$ be the field of fractions of $A$. By extension of scalars, 
we obtain a $K$-algebra $\cH_K=K\otimes_A \cH$. This algebra is known to 
be split semisimple; see \cite[9.3.5]{GePf}. (Note the form of the quadratic
relations for $T_s$, $s \in S$.) Furthermore, by Tits' Deformation Theorem, 
the irreducible representations of $\cH_K$ (up to isomorphism) are in 
bijection with the irreducible representations of $W$; see \cite[8.1.7]{GePf}. 
Thus, we can write
\[ \Irr(\cH_K)=\{E^\lambda_v \mid \lambda \in \Lambda\}.\]
The correspondence $E^\lambda \leftrightarrow E^\lambda_v$ is uniquely 
determined by the following condition:
\[ \mbox{trace}\bigl(w,E^\lambda\bigr)=\mbox{trace}\bigl(
T_w,E^\lambda_v\bigr)\Big|_{v=1} \qquad \mbox{for all $w \in W$};\]
note that $\mbox{trace}\bigl(T_w,E^\lambda_v\bigr) \in A$ for all 
$w\in W$. See also \cite[20.2, 20.3]{Lusztig03} for a discussion of the above
correspondence, but note that this relies on the validity of Lusztig's 
conjectures (see \S \ref{sub24} below).

\subsection{The integers $\ba_\lambda$ and $f_\lambda$} \label{sub22} 

The algebra $\cH$ is {\em symmetric}, where $\{T_w \mid w \in W\}$ and
$\{T_{w^{-1}}\mid w \in W\}$ form a pair of dual bases. Hence we
have the following orthogonality relations for the irreducible
representations of $\cH_K$:
\[ \sum_{w \in W} \mbox{trace}\bigl(T_w,E^\lambda_v\bigr)
\,\mbox{trace}\bigl(T_{w^{-1}},E_v^\mu\bigr)=\left\{\begin{array}{cl}
d_\lambda\,\bc_\lambda & \quad \mbox{if $\lambda=\mu$},\\ 0 & \quad 
\mbox{if $\lambda \neq\mu$};\end{array}\right.\]
see \cite[8.1.7]{GePf}. Here, $0 \neq \bc_\lambda \in {\Z}[v,v^{-1}]$ and,
following Lusztig, we can write
\[ \bc_\lambda=f_\lambda\, v^{-2\ba_\lambda}+\mbox{combination of strictly
higher powers of $v$},\]
where $\ba_\lambda,f_\lambda$ are integers such that $\ba_\lambda\geq 0$ 
and $f_\lambda>0$; see \cite[9.4.7]{GePf}. Thus, using $\cH$, we have 
associated with each $E^\lambda \in \Irr(W)$ two integers $\ba_\lambda$ 
and $f_\lambda$. These integers are explicitly known for all types of $W$; 
see Lusztig \cite[Chap.~22]{Lusztig03}.

Now let $p$ be a prime number. We say that $p$ is $L$-{\em bad} for $W$
if $p$ divides $f_\lambda$ for some $\lambda \in \Lambda$. Otherwise,
$p$ is called $L$-{\em good}. If $L$ is a positive multiple of the length 
function, this corresponds to the familiar definition of ``bad'' primes; 
see Lusztig \cite[Chap.~4]{Lu4}. Recall that, in this case, the conditions 
for being good the various irreducible types of $W$ are as follows:
\begin{center} $\begin{array}{rl} A_n: & \mbox{no condition}, \\
B_n, C_n, D_n: & p \neq 2, \\
G_2, F_4, E_6, E_7: &  p \neq 2,3, \\
E_8: & p \neq 2,3,5.  \end{array}$
\end{center}
See \cite[Example~4.7]{mylaus} and \cite[Def.~2.3]{GeJa} for a description 
of the $L$-bad primes when $L$ is not constant on $S$.

\subsection{The Kazhdan--Lusztig basis of $\cH$} \label{sub23} 

Let $\{c_w\mid w\in W\}$ be the ``new'' basis of $\cH$ defined in 
\cite[Theorem~5.2]{Lusztig03}. We have $c_w=T_w+ \sum_{y} p_{y,w}\, T_y$ 
where $p_{y,w} \in v^{-1}{\Z}[v^{-1}]$ and $p_{y,w}=0$ unless $y<w$ in 
the Bruhat--Chevalley order. Given $x,y\in W$, we write
\[c_xc_y=\sum_{z\in W} h_{x,y,z}c_z\qquad\mbox{where $h_{x,y,z}\in A$}.\]
As in \cite{Lusztig03}, we usually work with the elements $c_w^\dagger$ 
obtained by applying the unique $A$-algebra involution $H \rightarrow H$,
$h \mapsto h^\dagger$ such that $T_s^\dagger=-T_s^{-1}$ for any $s\in S$;
see \cite[3.5]{Lusztig03}. We refer to \cite[Chap.~8]{Lusztig03} for the 
definition of the preorders $\leq_{\cL}$, $\leq_{\cR}$, $\leq_{\cL\cR}$ 
and the corresponding equivalence relations $\sim_{\cL}$, $\sim_{\cR}$, 
$\sim_{\cL\cR}$ on $W$, induced by $L$. The equivalence classes with 
respect to these relations are called {\em left}, {\em right} and 
{\em two-sided cells} of $W$, respectively. If $L(s)=1$ for all $s \in S$,
then $C_w=(-1)^{l(w)}c_w^\dagger$ where $\{C_w\}$ is the basis originally
constructed in \cite{KaLu}.

Let $z \in W$. Following Lusztig \cite[13.6]{Lusztig03}, we define 
$\ba(z)\in {\Z}_{\geq 0}$ by the condition that 
\begin{alignat*}{2}
v^{\ba(z)}h_{x,y,z} &\in {\Z}[v] &&\qquad \mbox{for all $x,y\in W$},\\
v^{\ba(z)-1}h_{x,y,z} &\not\in {\Z}[v] &&\qquad \mbox{for some 
$x,y\in W$}.
\end{alignat*}
Furthermore, if $p_{1,z}\neq 0$, we define $\Delta(z)\in \Z_{\geq 0}$ 
and $0\neq n_z\in \Z$ by the condition that
\[ p_{1,z}=n_zv^{-\Delta(z)}+\mbox{strictly smaller  powers of $v$};
\quad \mbox{see \cite[14.1]{Lusztig03}}.\]
Otherwise, we set $\Delta(z)=\infty$ and leave $n_z$ undefined. (This case
can only occur when $L(s)=0$ for some $s \in S$.) We set 
\[ \cD=\{z \in W \mid \ba(z)=\Delta(z)\}.\]
(In what follows, the coefficients $n_z$ will only play a role when 
$z\in \cD$.)

\subsection{Lusztig's conjectures and the asymptotic ring $J$} \label{sub24}

In the sequel, we assume that the following hypotheses are satisfied.
\begin{center}
\fbox{\em Lusztig's conjectures {\bf (P1)}--{\bf (P15)} in 
\cite[14.2]{Lusztig03} hold for $\cH$.}
\end{center}
(Actually, for our purposes, instead of {\bf (P15)} it is enough to 
require the somewhat weaker statement \cite[18.9(b)]{Lusztig03}, which 
we called {\bf (P15')} in \cite[\S 5]{mylaus}.)
By \cite[Chap.~15]{Lusztig03}, these conjectures are known to hold when $L$ 
is a positive multiple of the length function (``equal parameter case''), 
thanks to a deep geometric interpretation of the basis $\{c_w\}$. They are 
also known for a certain class of non-trivial weight functions in type $B_n$; 
see \cite{GeIa06}, \cite{myert06}.

Assuming the above hypotheses, we can perform the following constructions.
Following Lusztig \cite[Chap.~18]{Lusztig03}, let $J$ be a free 
${\Z}$-module with basis $\{t_w\mid w\in W\}$. We define a multiplication on 
$J$ by
\[ t_xt_y=\sum_{z\in W} \gamma_{x,y,z^{-1}}\, t_z,\]
where $\gamma_{x,y,z^{-1}} \in \Z$ is the constant term of $v^{\ba(z)}
h_{x,y,z} \in {\Z}[v]$. Then it turns out that this multiplication is 
associative and we have an identity element given by $1_J=\sum_{d\in \cD} 
n_dt_d$. Furthermore, we have a homomorphism of $A$-algebras $\phi \colon 
\cH \rightarrow J_A=A \otimes_{\Z} J$ such that 
\[ \phi(c_w^\dagger)=\sum_{\atop{z\in W,d\in \cD}{\ba(z)=\ba(d)}}
h_{w,d,z} \hat{n}_z\, t_z \qquad (w\in W),\]
where $\hat{n}_z$ is defined as follows. Given $z\in W$, let 
$d$ be the unique element of $\cD$ such that $d\sim_{\cL} z^{-1}$; then
$\hat{n}_z=n_d=\pm 1$. (See properties {\bf (P5)}, {\bf (P13)} in 
Lusztig's conjectures \cite[14.2]{Lusztig03}.)  Note that the function 
$z \mapsto \hat{n}_z$ is constant on the right cells of $W$.

\subsection{Representations of $J$} \label{sub25} 

Until the end of this section, we will assume that $R=\Q$. Upon 
substituting $v\mapsto 1 \in \Q$, the algebra $\cH$ specialises to 
${\Q}[W]$. Hence, we obtain a homomorphism of $\Q$-algebras
\[ \phi_1\colon {\Q}[W]  \rightarrow J_{\Q}, \qquad \mbox{where $J_{\Q}=
{\Q} \otimes_{\Z} J$}.\]
This is an isomorphism by the argument in \cite[20.1]{Lusztig03}. Since 
${\Q}[W]$ is split semisimple,  we can conclude that $J_{\Q}$ also is 
split semisimple. Via the isomorphism $\phi_1 \colon {\Q}[W]\rightarrow 
J_\Q$, the set $\Lambda$ can be used to parametrize the irreducible 
representations of $J_\Q$ (up to isomorphism). As in 
\cite[20.2]{Lusztig03}, we write
\[ \Irr(J_{\Q})=\{E_\spadesuit^\lambda \mid \lambda \in \Lambda\},\]
where $E_\spadesuit^\lambda$ coincides with $E^\lambda$ as $\Q$-vector space 
and the action of $j \in J_\Q$ on $E_\spadesuit^\lambda$ is the same as the 
action of $\phi_1^{-1}(j)$ on $E^\lambda$. Now, by \cite[20.1(b)]{Lusztig03}, 
the algebra $J$ is symmetric where $\{t_w \mid w \in W\}$ and $\{t_{w^{-1}} 
\mid w \in W\}$ form a pair of dual bases. Thus, we have the following
orthogonality relations:
\[ \sum_{w \in W} \mbox{trace}\bigl(t_w,E_\spadesuit^\lambda\bigr)
\,\mbox{trace}\bigl(t_{w^{-1}},E_\spadesuit^\mu\bigr)= \left\{
\begin{array}{cl} d_\lambda\,f_\lambda & \quad \mbox{if $\lambda=\mu$},\\ 
0 & \quad \mbox{if $\lambda \neq\mu$}.\end{array}\right.\]
(The fact that $f_\lambda$ appears on the right hand side is shown in 
\cite[20.11]{Lusztig03}.) Choosing a vector space basis of $E^\lambda$, 
we obtain a matrix representation 
\[ \rho^\lambda \colon J_{\Q} \rightarrow M_{d_\lambda}(\Q), \qquad
t_w \mapsto \bigl(\rho_{\fs\ft}^\lambda(t_w)\bigr)_{\fs,\ft\in M(\lambda)},\]
where $M(\lambda)=\{1,\ldots,d_\lambda\}$. Then the following {\em Schur 
relations} hold; see \cite[\S 7.2]{GePf}. For $\lambda, \mu \in 
\Lambda$, $\fs,\ft \in M(\lambda)$ and $\fu,\fv\in M(\mu)$, we have:
\[ \sum_{w \in W} \rho_{\fs\ft}^\lambda(t_{w})\,\rho^\mu_{\fu\fv}(
t_{w^{-1}})= \left\{\begin{array}{cl} f_\lambda & \quad \mbox{if 
$(\lambda,\fs,\ft)= (\mu,\fv,\fu)$},\\ 0 & \quad \mbox{otherwise}.
\end{array}\right.\]
These equations can be inverted and this yields the ``second'' Schur 
relations:
\[ \sum_{\lambda \in \Lambda} \sum_{\fs,\ft \in M(\lambda)}
\frac{1}{f_\lambda}\rho_{\fs\ft}^\lambda(t_x)\,\rho_{\ft\fs}^\lambda(
t_{y^{-1}})= \left\{\begin{array}{cl} 1 & \quad \mbox{if $x=y$},\\ 0 & 
\quad \mbox{otherwise}.\end{array}\right.\]
By a general argument, every irreducible representation of $J_\Q$ 
leaves a positive-definite quadratic form invariant. More precisely, we 
have the following result.

\addtocounter{theorem}{5}
\begin{proposition} \label{lem1} Let $\lambda \in \Lambda$. Then a basis of
$E^\lambda$ can be chosen such that
\begin{itemize}
\item[(a)] $\rho^\lambda_{\fs\ft}(t_w) \in \Z$ for all $w\in W$
and $\fs,\ft \in M(\lambda)$.
\end{itemize}
Furthermore, there exists a symmetric, positive-definite matrix 
\[B^\lambda=(\beta_{\fs\ft}^\lambda)_{\fs,\ft\in M(\lambda)} \qquad 
\mbox{where} \qquad\beta_{\fs\ft}^\lambda\in \Z \mbox{ for all $\fs,\ft 
\in M_{d_\lambda}$},\]
such that the following two conditions hold:
\begin{itemize}
\item[(b)] $B^\lambda\cdot \rho^\lambda(t_{w^{-1}})=\rho^\lambda
(t_w)^{\operatorname{tr}}\cdot B^\lambda$ for all $w \in W$;  
\item[(c)] Any prime which divides $\det(B^\lambda)\neq 0$ is $L$-bad 
for $W$.
\end{itemize}
\end{proposition}

(For any matrix $M$, we denote by $M^{\text{tr}}$ the 
transpose of $M$.)

\begin{proof}
To simplify the notation, write $E=E^\lambda$, $d=d_\lambda$ and  
$\rho=\rho^\lambda$ so that we can omit a subscript or superscript 
$\lambda$ in the subsequent formulas. Since $J$ is defined over $\Z$, the 
statement in (a) follows by a general argument; see, for example, 
\cite[7.3.7]{GePf}.  Now let 
\[B_1=\sum_{y\in W} \rho(t_y)^{\operatorname{tr}}\cdot\rho(t_y)\in 
M_d(\Z).\]
This matrix clearly is symmetric.  Now let $0\neq e=(e_1,\ldots,e_d)
\in \Z^d$. Since $\rho$ is irreducible, there exists some $y \in W$
such that $e\rho(t_y)^{\operatorname{tr}}\neq 0$ and, hence, the 
standard scalar product of this vector with itself will be strictly 
positive. Consequently, we have $eB_1e^{\operatorname{tr}}>0$. Thus, 
$B_1$ is positive-definite and, in particular, $\det(B_1)\neq 0$. For 
any $x\in W$, we have
\[ B_1\cdot \rho(t_{x^{-1}})=\sum_{y \in W} \rho(t_y)^{\operatorname{tr}}
\cdot \rho(t_yt_{x^{-1}})= \sum_{y,z\in W} \gamma_{y,x^{-1},z^{-1}}\,
\rho(t_y)^{\operatorname{tr}}\cdot\rho(t_z).\]
Now $\gamma_{y,x^{-1},z^{-1}}=\gamma_{x^{-1},z^{-1},y}$ by {\bf (P7)} in 
\cite[14.2]{Lusztig03}, and $\gamma_{x^{-1},z^{-1},y}=\gamma_{z,x,y^{-1}}$ 
by \cite[Prop.~13.9]{Lusztig03}. Hence, the right hand side of the 
above identity equals
\[\sum_{y,z\in W} \gamma_{z,x,y^{-1}} \rho(t_y)^{\operatorname{tr}} 
\cdot \rho(t_z)=\sum_{z\in W} \rho(t_zt_x)^{\operatorname{tr}} 
\cdot \rho(t_z)= \rho(t_x)^{\operatorname{tr}}\cdot B_1.\]
Let $0 \neq n \in \Z$ be the greatest common divisor of all non-zero
coefficients of $B_1$ and $B=n^{-1}B_1\in M_d(\Z)$. Then
$B$ is a symmetric, positive-definite matrix such that (b) holds. It
remains to prove (c). Let $p$ be a prime number and denote by $\bar{B}$ 
the matrix obtained by reducing all coefficents modulo $p$. By reduction 
modulo $p$, we also obtain an $\F_p$-algebra $J_p=\F_p \otimes_\Z J$ and 
a corresponding matrix representation $\bar{\rho}\colon J_p \rightarrow 
M_{d}(\F_p)$. Hence we have
\[ \bar{B}\neq 0 \qquad \mbox{and}\qquad \bar{B}\cdot\bar{\rho}(t_w)=
\bar{\rho}(t_{w^{-1}})^{\text{tr}}\cdot \bar{B}\quad\mbox{for all $w\in W$}.\]
By \cite[13.9]{Lusztig03}, the map $t_w\mapsto t_{w^{-1}}$ defines an 
involutory anti-auto\-morphism of $J$. Hence the assignment $t_w \mapsto 
\bar{\rho}(t_{w^{-1}})^{\text{tr}}$ also defines a representation of $J_p$. 
The above identity now shows that $\bar{B}\neq 0$ is an ``intertwining 
operator''. Hence, if we knew that $\bar{\rho}$ was irreducible, then Schur's 
Lemma would imply that $\bar{B}$ were invertible and so $p$ could not 
divide $\det(B)$. 

Thus, it remains to show that $\bar{\rho}$ is an irreducible representation 
of $J_p$ whenever $p$ is $L$-good. But this follows from a general 
argument about symmetric algebras. Indeed, as already noted in \S \ref{sub25}, 
$J$ is symmetric and we have the Schur relations for the matrix 
coefficients of $\rho=\rho^\lambda$. Reducing these relations modulo $p$,
we obtain:
\[ \sum_{w \in W} \bar{\rho}_{\fs\ft}^\lambda(t_w) \, 
\bar{\rho}^\lambda_{\fu\fv} (t_{w^{-1}})= \left\{\begin{array}{cl} 
f_\lambda\bmod p & \quad \mbox{if $\fs=\fv$, $\fu=\ft$},\\ 0 & \quad 
\mbox{otherwise}.\end{array}\right.\]
Since $f_\lambda \not \equiv 0 \bmod p$, one easily deduces from this that
$\bar{\rho}^\lambda$ is (absolutely) irreducible; see 
\cite[Remark~7.2.3]{GePf}. %\qed
\end{proof}

%%%%%%%%%%%%%%%%%%%%%%%%%%%%%%%%%%%%%%%%%%%%%%%%%%%%%%%%%%%%%%%%%%%%%%%%%%%%%
\section{A cell datum for $\cH$} \label{secDatum}
We keep the notation of the previous sections. In order to show that
$\cH$ is ``cellular'' in the sense of Graham--Lehrer 
\cite[Definition~1.1]{GrLe}, we must specify a quadruple $(\Lambda,M,C,*)$ 
satisfying the following conditions.
\begin{itemize}
\item[(C1)] $\Lambda$ is a partially ordered set, $\{M(\lambda) \mid
\lambda \in \Lambda\}$ is a collection of finite sets  and
\[ C \colon \coprod_{\lambda \in \Lambda} M(\lambda)\times M(\lambda)
\rightarrow \cH \]
is an injective map whose image is an $A$-basis of $\cH$;
\item[(C2)] If $\lambda \in \Lambda$ and $\fs,\ft\in M(\lambda)$, write
$C(\fs,\ft)=C_{\fs,\ft}^\lambda \in \cH$. Then $* \colon \cH \rightarrow 
\cH$ is an $A$-linear anti-involution such that $(C_{\fs,\ft}^\lambda)^*=
C_{\ft,\fs}^\lambda$.
\item[(C3)] If $\lambda \in \Lambda$ and $\fs,\ft\in M(\lambda)$, then for 
any element $h \in \cH$ we have
\[ hC_{\fs,\ft}^\lambda\equiv \sum_{\fs'\in M(\lambda)} r_h(\fs',\fs)\,
C_{\fs',\ft}^\lambda\quad \bmod \cH(<\lambda),\]
where $r_h(\fs',\fs) \in A$ is independent of $\ft$ and where $\cH(<\lambda)$
is the $A$-submodule of $\cH_n$ generated by $\{C_{\fs'',\ft''}^\mu
\mid \mu <\lambda; \fs'',\ft''\in M(\mu)\}$.
\end{itemize}
We now define a required quadruple $(\Lambda,M,C,*)$ as follows.

Let $\Lambda$ be an indexing set for the irreducible representations of
$W$, as in \S \ref{sub21}. Using the $\ba$-invariants in \S \ref{sub22}, we 
define a partial order $\preceq$ on $\Lambda$ by
\[ \lambda \preceq \mu \qquad \stackrel{\text{def}}{\Leftrightarrow} \qquad
\lambda=\mu \quad \mbox{or}\quad \ba_\lambda >\ba_\mu.\]
Thus, $\Lambda$ is ordered according to {\em decreasing} $\ba$-value. Next, 
we define an $A$-linear anti-involution $* \colon \cH \rightarrow \cH$ 
by $T_w^*=T_{w^{-1}}$ for all $w \in W$. Thus, $T_w^*=T_w^{\,\flat}$ in the
notation of \cite[3.4]{Lusztig03}.

For $\lambda\in \Lambda$, we set $M(\lambda)=\{1,\ldots,d_\lambda\}$ as 
before. The trickiest part is, of course, the definition of the basis 
elements $C_{\fs,\ft}^\lambda$ for $\fs,\ft\in M(\lambda)$. We can now 
state the main result of this paper.

\begin{theorem} \label{mainthm} Recall that Lusztig's conjectures 
{\bf (P1)}--{\bf (P15)} are assumed to hold for $\cH$. Let $R \subseteq \C$ 
be a subring such that all $L$-bad primes are invertible in $R$. Let
$\bigl(\rho_{\fs\ft}^\lambda(t_w)\bigr)$ and  $\bigl(\beta_{\fs
\ft}^\lambda\bigr)$ be as in Proposition~\ref{lem1}. For any $\lambda 
\in \Lambda$ and $\fs,\ft\in M(\lambda)$, define
\[C_{\fs,\ft}^\lambda=\sum_{w \in W} \sum_{\fu \in M(\lambda)} \hat{n}_w
\hat{n}_{w^{-1}}\,\beta^\lambda_{\ft\fu}\,\rho_{\fu\fs}^\lambda\, 
(t_{w^{-1}})\, c_w^\dagger.\]
Then $C_{\fs,\ft}^\lambda$ is a $\Z$-linear combination of Kazhdan--Lusztig 
basis elements $c_w^\dagger$ where $\ba(w)=\ba_\lambda$. The quadruple 
$(\Lambda,M,C,*)$ is a ``cell datum'' in the sense of Graham--Lehrer 
\cite{GrLe}. 
\end{theorem}

\begin{proof} First note that, by \cite[Prop.~20.6]{Lusztig03}, we have
$\rho^\lambda(t_{w^{-1}})=0$ unless $\ba(w)=\ba_\lambda$. Thus, 
$C_{\fs,\ft}^\lambda$ is an integral linear combination of elements 
$c_w^\dagger$ where $\ba(w)=\ba_\lambda$. In what follows, it will
be convenient to write the coefficients occurring in various sums as
entries of matrices. Thus, for example, the defining formula for 
$C_{\fs,\ft}^\lambda$ reads:
\begin{align*}
C_{\fs,\ft}^\lambda&=\sum_{w \in W} \hat{n}_w\hat{n}_{w^{-1}}\,\bigl(
B^\lambda \cdot\rho^\lambda(t_{w^{-1}})\bigr)_{\ft,\fs}\, c_w^\dagger\\
&= \sum_{w \in W} \hat{n}_w\hat{n}_{w^{-1}}\,\bigl(\rho^\lambda
(t_w)^{\text{tr}}\cdot B^\lambda\bigr)_{\ft,\fs}\, c_w^\dagger,
\end{align*}
where the second equality holds by Proposition~\ref{lem1}. 
We now proceed in three steps.

\medskip
{\em Step 1}. (C1) holds, that is, the elements $\{C_{\fs,\ft}^\lambda 
\mid \lambda \in \Lambda \mbox{ and }  \fs,\ft\in M(\lambda)\}$ form a 
basis of $\cH$. This is proved as follows. By Wedderburn's Theorem, 
$\dim \cH_K= |W|=\sum_{\lambda \in \Lambda} |M(\lambda)|^2$.  Hence the 
above set has the correct cardinality. It is now sufficient to show that the 
elements $\{C_{\fs,\ft}^\lambda\}$ span $\cH$ as an $A$-module.

Let us fix $y \in W$.  We consider the following $R$-linear combination:
\[\sum_{\lambda \in \Lambda} \sum_{\fs,\ft\in M(\lambda)} \frac{1}{f_\lambda} 
\,\bigl(\rho^\lambda(t_y)\cdot (B^\lambda)^{-1}\bigr)_{\fs,\ft}\, 
C_{\fs,\ft}^\lambda.\]
Note that the coefficients lie in $R$ since $f_\lambda$ and $\det(B^\lambda)$
are invertible in $R$. Inserting the first of the above-mentioned two
expressions for $C_{\fs,\ft}^\lambda$, we obtain:
\[ \sum_{w\in W} \hat{n}_w\hat{n}_{w^{-1}}
\sum_{\lambda \in \Lambda} \sum_{\fs\in M(\lambda)} \frac{1}{f_\lambda} \,
\bigl(\rho^\lambda(t_y)\cdot \rho^\lambda(t_{w^{-1}})\bigr)_{\fs,\fs}\,
c_w^\dagger.\]
Now, writing out the product $\rho^\lambda(t_y)\cdot 
\rho^\lambda(t_{w^{-1}})$ and using the ``second'' Schur relations, 
we find that 
\[ \sum_{\lambda \in \Lambda} \sum_{\fs\in M(\lambda)} \frac{1}{f_\lambda} \,
\bigl(\rho^\lambda(t_y)\cdot \rho^\lambda(t_{w^{-1}})\bigr)_{\fs,\fs} =
\delta_{yw}.\]
Hence our linear combination reduces to $\hat{n}_y\hat{n}_{y^{-1}}\, 
c_y^\dagger=\pm c_y^\dagger$. Thus, $c_y^\dagger$ is an $R$-linear
combination of the elements $C_{\fs,\ft}^\lambda$, as required.
  
\medskip
{\em Step 2}. (C2) holds, that is,  we have 
$(C_{\fs,\ft}^\lambda)^*=C_{\ft,\fs}^\lambda$ for all $\lambda \in \Lambda$ 
and $\fs,\ft\in M(\lambda)$. 
This is seen as follows. By \cite[4.9 and 5.6]{Lusztig03}, we have 
$(c_w^\dagger)^*=(c_w^*)^\dagger= c_{w^{-1}}^\dagger$. Thus,  using
the above two expressions for $C_{\fs,\ft}^\lambda$, we obtain:
\begin{align*}
(C_{\fs,\ft}^\lambda)^* &= \sum_{w \in W} \hat{n}_w\hat{n}_{w^{-1}}\,\bigl(
B^\lambda \cdot \rho^\lambda(t_{w^{-1}})\bigr)_{\ft,\fs}\, 
c_{w^{-1}}^\dagger\\ &= \sum_{w \in W} \hat{n}_w\hat{n}_{w^{-1}}\,\bigl(
\rho^\lambda(t_w)^{\operatorname{tr}}\cdot B^\lambda\bigr)_{\ft,\fs}\, 
c_{w^{-1}}^\dagger\\ &= \sum_{w \in W} \hat{n}_w\hat{n}_{w^{-1}}\,
\bigl( B^\lambda \cdot \rho^\lambda(t_w)\bigr)_{\fs,\ft}\, 
c_{w^{-1}}^\dagger = C_{\ft,\fs}^\lambda,
\end{align*}
as required.
 
\medskip
{\em Step 3}. Finally, we consider the multiplication rule (C3).
By Lusztig \cite[18.10]{Lusztig03}, there is a natural left $J_A$-module 
structure on $\cH$ given by the formula
\[ t_x * c_w^\dagger=\sum_{z \in W} \gamma_{x,w,z^{-1}}\, \hat{n}_w\,
\hat{n}_z\, c_z^\dagger \qquad (x,w \in W).\]
We begin by studying the effect of the $J_A$-action on the element
$C_{\fs,\ft}^\lambda$. Let $x \in W$.  We claim that 
\[ t_x * C_{\fs,\ft}=\sum_{\fs'\in M(\lambda)} \rho^\lambda_{\fs'\fs}(t_x)
\, C_{\fs', \ft}.\]
Indeed, recalling the defining formula for $C_{\fs,\ft}^\lambda$, we have:
\begin{align*}
t_x * C_{\fs,\ft}^\lambda &=
\sum_{w \in W} \sum_{\fu \in M(\lambda)} \hat{n}_w \hat{n}_{w^{-1}}\, 
\beta_{\ft\fu}^\lambda\, \rho_{\fu\fs}^\lambda(t_{w^{-1}})\, t_x * 
c_w^\dagger\\ &=\sum_{w,z \in W} \sum_{\fu \in M(\lambda)} \hat{n}_z
\hat{n}_{w^{-1}}\, \beta_{\ft\fu}^\lambda\, \rho_{\fu\fs}^\lambda(
t_{w^{-1}}) \, \gamma_{x,w,z^{-1}} \, c_z^\dagger.
\end{align*}
Assume that the term in the sum corresponding to $w,z\in W$ is non-zero.
Then $\gamma_{x,w,z^{-1}}\neq 0$ and so $w \sim_{\cL} z$, by property 
{\bf (P8)} in \cite[14.2]{Lusztig03}. Hence $w^{-1}\sim_{\cR}z^{-1}$ and
so $\hat{n}_{w^{-1}}=\hat{n}_{z^{-1}}$. Furthermore, $\gamma_{x,w,z^{-1}}
=\gamma_{z^{-1},x,w}$ by property {\bf (P7)} in \cite[14.2]{Lusztig03}. 
Thus, we obtain 
\begin{align*}
t_x * C_{\fs,\ft}^\lambda &=
\sum_{z \in W} \sum_{\fu \in M(\lambda)} \hat{n}_z\hat{n}_{z^{-1}}\,
\beta_{\ft\fu}^\lambda\, \rho_{\fu\fs}^\lambda\Bigl(\sum_{w \in W}
\gamma_{z^{-1},x,w} t_{w^{-1}}\Bigr)\,c_z^\dagger\\
&=\sum_{z \in W} \sum_{\fu\in M(\lambda)}\hat{n}_z\hat{n}_{z^{-1}}\,
\beta_{\ft\fu}^\lambda\,\rho_{\fu\fs}^\lambda\bigl(t_{z^{-1}}\, t_x\bigr)
\,c_z^\dagger\\
&=\sum_{z \in W} \hat{n}_z\hat{n}_{z^{-1}}\,\bigl(B^\lambda\cdot
\rho^\lambda(t_{z^{-1}}t_x)\bigr)_{\ft,\fs}\, c_z^\dagger\\ 
&=\sum_{z \in W} \sum_{\fs'\in M(\lambda)}\hat{n}_z\hat{n}_{z^{-1}}\,
\bigl(B^\lambda\cdot \rho^\lambda(t_{z^{-1}})\bigr)_{\ft,\fs'}\,
\rho_{\fs'\fs}^\lambda(t_x)\, c_z^\dagger\\ &=\sum_{\fs'\in M(\lambda)} 
\rho^\lambda_{\fs'\fs}(t_x)\, C_{\fs',\ft},
\end{align*}
as claimed.
Now let $h \in \cH$ and write $\phi(h)=\sum_{x \in W} a_h(x)\,t_x$ where
$a_h(x) \in A$. Then define 
\[ r_h(\fs',\fs)=\sum_{x \in W} a_h(x)\, \rho^\lambda_{\fs'\fs}(t_x)
 \quad \mbox{for any $\fs',\fs \in M(\lambda)$}.\]
Note that, indeed, this coefficient lies in $A$ and it only depends on
$\fs,\fs'$ and $h$. Then the above computation shows that 
\[ \phi(h) * C_{\fs,\ft}^\lambda=\sum_{\fs' \in M(\lambda)} r_h(\fs',\fs)\, 
C_{\fs',\ft}^\lambda.\]
For any $a \geq 0$, we define $\cH^{\geq a}$ to be the $A$-span  of
all elements $c_y^\dagger$ where $y \in W$ is such that $\ba(y)\geq a$.
By \cite[18.10(a)]{Lusztig03}, we have
\[hc_w^\dagger\equiv\phi(h)*c_w^\dagger\quad\bmod\cH^{\geq \ba(w)+1}
\quad \mbox{for any $h \in \cH$ and $w \in W$}.\]
We have already noted in the beginning of the proof that 
$C_{\fs,\ft}^\lambda$ is a linear combination of elements $c_z^\dagger$ 
where $\ba(z)= \ba_\lambda$. Hence the above relations imply that
\[ hC_{\fs,\ft}^\lambda \equiv \phi(h)*C_{\fs,\ft}^\lambda \equiv 
\sum_{\fs'\in M(\lambda)} r_h(\fs',\fs)\, C_{\fs',\ft}^\lambda \quad \bmod 
\cH^{\geq \ba_{\lambda}+1}.\]
The definition of the partial order $\preceq$ on $\Lambda$ now shows
that (C3) holds. %\qed
\end{proof}

\begin{corollary} \label{cor1} Let $\theta \colon A \rightarrow k$ be a ring
homomorphism into an integral domain $k$ (i.e., a ``specialisation''). 
By extension of scalars, we obtain a $k$-algebra $\cH_k=k \otimes_A \cH$. 
Then the above ingredients define a ``cell datum'' for $\cH_k$, where
\[ C_{\fs,\ft}^\lambda=\sum_{w \in W} \sum_{\fu\in M(\lambda)} 
\hat{n}_w\hat{n}_{w^{-1}}\,\theta(\beta^\lambda_{\ft\fu})\,\theta
\bigl(\rho_{\fu\fs}^\lambda\, (t_{w^{-1}}) \bigr)\, (1 \otimes c_w^\dagger)
\in \cH_k,\]
for $\lambda \in \Lambda$ and $\fs,\ft \in M(\lambda)$.
\end{corollary}

\begin{proof} This immediately follows from Theorem~\ref{mainthm}; see the
remarks in \cite[(1.8)]{GrLe}. %\qed 
\end{proof}

Thus, since Lusztig's conjectures {\bf (P1)}--{\bf (P15)} hold
in the ``equal parameter'' case, we have proved Theorem~\ref{thmb},
as stated in the introduction.

\begin{remark} \label{noncryst} Assume that $W$ is a finite Coxeter group which 
is not a Weyl group, i.e., $W$ is of type $H_3$, $H_4$ or $I_2(m)$ where 
$m=5$ or $m\geq 7$. In these cases, $\Q$ is no longer a splitting field 
for $W$, and this leads to some technical complications. However, choosing 
$R\subseteq \C$ appropriately, Theorem~\ref{mainthm} remains valid in these 
cases as well. 

First note that Lusztig's properties {\bf (P15)}--{\bf (P15)} (see 
\S \ref{sub24}) are known to hold for types $H_3$, $H_4$ and $I_2(m)$ 
(where $L$ is a positive multiple of the length function), thanks to 
DuCloux \cite{Fokko}. (As far as we know, these properties do not seem to 
have been verified in type $I_2(m)$ with unequal parameters.) 

Now, the irreducible representations in type $H_3$, $H_4$ and $I_2(m)$ (and 
much further information) are explicitly known; see 
\cite[\S 8.1 and \S 11.2]{GePf}. Using this explicit information, one can 
show that Theorem~\ref{mainthm} holds if we require that the subring 
$R \subseteq \C$ satisfies the following two conditions: (1) $R$ contains 
$\zeta+\zeta^{-1}$ where $\zeta$ is a root of unity of order $5$
(in type $H_3$, $H_4$) or $m$ (in type $I_2(m)$); (2) the order of $W$
is invertible in $R$.  

A more detailed discussion of these cases will appear elsewhere. 

Jeong et al.\  \cite{kor} point out that Fakiolas \cite{Fak}  also produces 
a cellular basis in type $I_2(m)$ where $L$ is a positive multiple of the 
length function. 
\end{remark}

%%%%%%%%%%%%%%%%%%%%%%%%%%%%%%%%%%%%%%%%%%%%%%%%%%%%%%%%%%%%%%%%%%%%%%%%%%%%%
\section{Examples and applications to modular representations} \label{secApp}

Throughout this section, we assume that Lusztig's conjectures on Hecke
algebras with unequal parameters hold; see \S \ref{sub24}. (Recall that
this is the case, for example, if the weight function $L$ is a positive
multiple of the length function.) We now discuss examples and applications 
of Theorem~\ref{mainthm}.

\begin{example} \label{expdim1} Assume that $L(s)>0$ for all $s \in S$ and
that $\lambda \in \Lambda$ is such that $\dim E^\lambda=d_\lambda=1$. Then 
there exists a group homomorphism $\eta \colon W \rightarrow\{\pm 1\}$ and 
a weight function $m \colon W \rightarrow \Z$ such that $T_w$ acts on 
$E_v^\lambda$ via the $A$-algebra homomorphism 
\[\mu^\lambda\colon\cH\rightarrow A,\qquad T_w\mapsto\eta(w)v^{m(w)}.\]
By \cite[Prop.~20.6]{Lusztig03}, the corresponding representation of 
$J$ is given by
\[ t_w \mapsto \left\{\begin{array}{cl}   (-1)^{l(w)}\,\hat{n}_w
\hat{n}_{w^{-1}}\, \eta(w) & \qquad \mbox{if $\ba_\lambda+m(w)=0$},\\ 
0 & \qquad \mbox{otherwise}.  \end{array}\right.\]
(Note that a factor $\hat{n}_w\hat{n}_{w^{-1}}$ should be inserted into
the formula in \cite[20.6]{Lusztig03}.) Clearly, we can take $B^\lambda=(1)$
in this case. Hence, the unique element of the cellular basis corresponding 
to $\lambda$ is given by 
\[ C_{1,1}^\lambda=\sum_{\atop{w \in W}{\ba_\lambda+m(w)=0}} 
(-1)^{l(w)}\,\eta(w)\, c_w^\dagger.\] 
For example, if $E^\lambda$ is the unit representation, we have
$\eta(w)=1$ and $m(w)=L(w)$ for all $w \in W$; furthermore, $\ba_\lambda=0$ 
and so $C_{1,1}^\lambda=c_1^\dagger$ (since $L(w)>0$ for all $w\neq 1$). 

If $E^\lambda$ is the sign representation, we have $\eta(w)=(-1)^{l(w)}$ 
and $m(w)=-L(w)$ for all $w \in W$. Furthermore, $\ba_\lambda=L(w_0)$ 
where $w_0\in W$ is the longest element; see \cite[20.18]{Lusztig03}. 
Since $L(w)<L(w_0)$ for all $w\neq w_0$, we obtain $C_{1,1}^\lambda=
c_{w_0}^\dagger$ in this case.
\end{example}

\begin{example} \label{expirred} Assume that there are no $L$-bad primes
for $W$, that is, we have $f_\lambda=1$ for all $\lambda \in \Lambda$. 
Then we claim that there are signs $\delta_w=\pm 1$ ($w\in W$) such that
\[ \{\delta_w\, c_w^\dagger\mid w \in W\} \qquad \mbox{is a cellular basis}.\]
This is seen as follows. Let $\Gamma$ be a left cell of $W$ and 
$[\Gamma]=\langle t_x \mid x \in \Gamma\rangle_\Q \subseteq J_{\Q}$. 
Then $[\Gamma]$ is a simple left ideal of $J_\Q$; see 
\cite[Chap.~21]{Lusztig03} and \cite[Cor.~4.8]{myert05}. Hence $[\Gamma]$ 
affords an irreducible representation of $J_\Q$, and all irreducible 
representations arise in this way. So, for any $\lambda \in \Lambda$, we 
can choose a left cell $\Gamma^\lambda$ such that $[\Gamma^\lambda]\cong 
E_\spadesuit^\lambda$. Let us write 
\[ \Gamma^\lambda=\{x_\fs \mid \fs\in M(\lambda)\} \qquad \mbox{for
$\lambda \in \Lambda$}.\]
Then $\{t_{x_\fs} \mid \fs \in M(\lambda)\}$ is a basis of 
$[\Gamma^\lambda]$ and the corresponding matrix coefficients are
\[\rho^\lambda_{\fs\ft}(t_w)=\gamma_{w,x_\ft,x_{\fs}^{-1}}
\qquad \mbox{for $w \in W$ and $\fs,\ft\in M(\lambda)$}.\]
By {\bf (P7)} in \cite[14.2]{Lusztig03} and \cite[Prop.~13.9]{Lusztig03}, 
we have 
\[ \gamma_{w,x_\ft,x_{\fs}^{-1}}=\gamma_{x_{\ft}^{-1},w^{-1},x_\fs}=
\gamma_{w^{-1}, x_\fs,x_{\ft}^{-1}}.\]
This implies that $\rho^\lambda (t_{w^{-1}})=\rho^\lambda(t_w)^{\text{tr}}$
for all $w \in W$. Thus, if we take for $B^\lambda$ the identity matrix 
of size $d_\lambda$, then the conditions in Proposition~\ref{lem1} are
satisfied. Furthermore, for fixed $\fs,\ft\in M(\lambda)$, the Schur 
relations now read:
\[\sum_{y \in W} \bigl(\gamma_{y,x_\ft,x_\fs^{-1}}\bigr)^2=\sum_{y \in W} 
\rho_{\fs\ft}^\lambda(t_y)\,\rho_{\ft\fs}^\lambda(t_{y^{-1}})=f_\lambda=1.\]
We deduce that there is a unique $w=w_\lambda(\fs,\ft)\in W$ such that
$\gamma_{w,x_\ft,x_\fs^{-1}}\neq 0$; in fact, we have $\gamma_{w,x_\ft,
x_\fs^{-1}}= \pm 1$. Now the formula in Theorem~\ref{mainthm} reads 
\begin{align*}
C_{\fs,\ft}^\lambda&=\sum_{y \in W} \hat{n}_y\hat{n}_{y^{-1}}\,
\gamma_{y^{-1},x_\fs,x_\ft^{-1}}\, c_y^\dagger= \sum_{y \in W} 
\hat{n}_y\hat{n}_{y^{-1}}\,\gamma_{y,x_\ft,x_\fs^{-1}}\,
c_y^\dagger\\ &=\hat{n}_w\hat{n}_{w^{-1}}\,\gamma_{w,x_\ft,x_\fs^{-1}}\, 
c_w^\dagger=\delta_w\, c_w^\dagger,
\end{align*}
where $w=w_\lambda(\fs,\ft)$ and $\delta_w=\pm 1$. Thus, 
$\{\delta_w\,c_w^\dagger \mid w \in W\}$ is a cellular basis, as claimed. 

The above assumptions are satisfied for $W$ of type $A_n$ where $L$ is a
positive multiple of the length function; see \cite[22.4]{Lusztig03}. 
Thanks to a geometric interpretation of the basis $\{c_w\}$, we know that
$h_{x,y,z}\in {\Z}_{\geq 0}[v,v^{-1}]$ for all $x,y,z\in W$; see
\cite[Chap.~15]{Lusztig03} and the references there. Hence we also
have $\hat{n}_x\geq 0$ and $\gamma_{x,y,z}\geq 0$ for all $x,y,z\in W$. So, 
in this case, $\{c_w^\dagger\mid w\in W\}$ is a cellular basis of $\cH$, 
as originally pointed out by Graham--Lehrer \cite[Example~1.2]{GrLe}.
Note that another cellular structure in type $A_n$ was constructed
by Murphy \cite{Mu1}, \cite{Mu2}, using  purely combinatorial methods.
The exact relation between the two cellular structures is determined
in \cite{mymurph}.
\end{example}

\begin{example} \label{expb2} Let $W$ be the Weyl group of type $B_2$
where $S=\{s_1,s_2\}$ and $(s_1s_2)^4=1$. Let $L\colon W \rightarrow\Z$
be any weight function such that $L(s_i)>0$ for $i=1,2$. If $L(s_1)\neq
L(s_2)$, then there are no $L$-bad primes for $W$. Furthermore, by 
\cite{GeIa06}, \cite{myert06},  Lusztig's  conjectures hold in this case.
(More generally, they hold for a certain class of weight functions in 
type $B_n$ for any $n \geq 2$.) Hence, in this case, $\{\delta_w
c_w^\dagger \mid w \in W\}$ is a cellular basis by Example~\ref{expirred}. 
By \cite[Cor.~6.4]{myert06}, we actually have $\delta_w=1$ for all $w\in W$.

Now assume that $L(s_1)=L(s_2)=1$ (``equal parameter case''). Then $2$ is 
the only bad prime. We have $\Irr(W)=\{{\bf 1},\varepsilon_1,\varepsilon_2,
\varepsilon,r\}$ where ${\bf 1}$ is the unit representation, $\varepsilon$
is the sign representation,  $\varepsilon_1$, $\varepsilon_2$ have dimension
one, and $r$ has dimension two. The $\ba$-invariants are $\ba_{\bf 1}=0$, 
$\ba_{\varepsilon_1}=\ba_{\varepsilon_2}=\ba_r=1$ and $\ba_{\varepsilon}=4$. 
A cellular basis as in Theorem~\ref{mainthm} is given as follows:
\begin{alignat*}{2}
C_{1,1}^{\bf 1}&=c_1^\dagger, \qquad\qquad &
C_{1,1}^r &= c_{s_1}^\dagger+c_{s_1s_2s_1}^\dagger,\\
C_{1,1}^\varepsilon&=c_{w_0}^\dagger, \qquad\qquad &
C_{1,2}^r &= -2c_{s_1s_2}^\dagger, \\ C_{1,1}^{\varepsilon_1}&=
c_{s_2}^\dagger-c_{s_2s_1s_2}^\dagger,\qquad\qquad &C_{2,1}^r &= 
-2c_{s_2s_1}^\dagger, \\ C_{1,1}^{\varepsilon_2}&=c_{s_1}^\dagger-
c_{s_1s_2s_1}^\dagger,\qquad\qquad &C_{2,2}^r &= 
2c_{s_2}^\dagger+2c_{s_2s_1s_2}^\dagger.
\end{alignat*}
Indeed, by Example~\ref{expdim1}, we already know that $c_1^\dagger$ and 
$c_{w_0}^\dagger$ belong to the cellular basis. Now consider the two
further one-dimensional representations. We fix the notation such that
$\varepsilon_1(T_{s_1})=\varepsilon_2(T_{s_2})=v$ and $\varepsilon_1
(T_{s_1})=\varepsilon_2(T_{s_2})=-v^{-1}$. Using the formula in 
Example~\ref{expdim1}, we obtain the above expressions for 
$C_{1,1}^{\varepsilon_1}$ and $C_{1,1}^{\varepsilon_2}$. Finally, consider 
the two-dimensional representation $r$. We have to determine a corresponding 
representation of $J$ such that the conditions in Proposition~\ref{lem1} 
hold. Now, the basis elements $\{c_w\}$ are explicitly determined in
\cite[Prop.~7.3]{Lusztig03}. From this, one easily deduces the left cells. 
They are given by 
\[\{1\},\quad\{s_1,s_2s_1,s_1s_2s_1\},\quad\{s_2,s_1s_2,s_2s_1s_2\},
\quad\{w_0=s_1s_2s_1s_2\};\]
see \cite[8.7]{Lusztig03}. Furthermore, one finds the following relations 
in $J$. We have $t_1^2=t_1$ and $t_1t_x=0$ for $x\neq 1$; we have $t_{w_0}^2
=t_{w_0}$ and $t_{w_0}t_x=0$ for $x \neq w_0$. Furthermore,
\begin{gather*}
t_{s_1}^2=t_{s_1},\qquad t_{s_1}t_{s_1s_2}=t_{s_1s_2},\qquad 
t_{s_1}t_{s_1s_2s_1}=t_{s_1s_2s_1},\\
t_{s_2}^2=t_{s_2},\qquad t_{s_2}t_{s_2s_1}=t_{s_2s_1},\qquad t_{s_2}
t_{s_2s_1s_2}= t_{s_2s_1s_2},\\ 
t_{s_1s_2}t_{s_2s_1}=t_{s_1}+ t_{s_1s_2s_1},\qquad t_{s_1s_2}
t_{s_2s_1s_2}=t_{s_1s_2},\\
t_{s_2s_1}t_{s_1s_2}=t_{s_2}+ t_{s_2s_1s_2},\qquad 
t_{s_2s_1}t_{s_1s_2s_1}=t_{s_2s_1},\\
t_{s_1s_2s_1}t_{s_1s_2s_1}=t_{s_1},\qquad t_{s_2s_1s_2}t_{s_2s_1s_2}=t_{s_2}.
\end{gather*}
All other products $t_xt_y$ are either zero or can be deduced from the
above list using the formulas $\gamma_{x,y,z}=\gamma_{y^{-1},x^{-1},z^{-1}}$
(see \cite[Prop.~13.9]{Lusztig03}) and $\gamma_{x,y,z}=\gamma_{y,z,x}$
(see {\bf (P7)} in \cite[14.2]{Lusztig03}). Note also that, by {\bf (P8)}
in \cite[14.2]{Lusztig03}, we have $\gamma_{x,y,z}=0$ unless $x \sim_{\cL}
y^{-1}$, $y \sim_{\cL} z^{-1}$ and $z \sim_{\cL} x^{-1}$. Now one readily 
checks that the following assignments define an irreducible two-dimensional 
representation of $J_\Q$:
\begin{gather*} 
t_{1} \mapsto \begin{bmatrix} 0 & 0 \\ 0 & 0 \end{bmatrix},\quad
t_{s_1} \mapsto \begin{bmatrix} 1 & 0 \\ 0 & 0 \end{bmatrix},\quad
t_{s_2s_1} \mapsto \begin{bmatrix} 0 & 0 \\ -1 & 0 \end{bmatrix},\quad
t_{s_1s_2s_1} \mapsto \begin{bmatrix} 1 & 0 \\ 0 & 0 \end{bmatrix},\\
t_{s_2} \mapsto \begin{bmatrix} 0 & 0 \\ 0 & 1 \end{bmatrix},\quad
t_{s_1s_2} \mapsto \begin{bmatrix} 0 & -2 \\ 0 & 0 \end{bmatrix},\quad
t_{s_2s_1s_2} \mapsto \begin{bmatrix} 0 & 0 \\ 0 & 1 \end{bmatrix},\quad
t_{w_0} \mapsto \begin{bmatrix} 0 & 0 \\ 0 & 0 \end{bmatrix}.
\end{gather*}
Finally, taking for $B^\lambda$ the diagonal matrix with diagonal entries
$1$ and $2$, we see that all the conditions in Proposition~\ref{lem1} 
are satisfied. This yields the expressions for the basis elements 
$C_{\fs,\ft}^r$ where $1 \leq \fs,\ft \leq 2$.
\end{example}
%\begin{alignat*}{2}
%{\bf 1}\colon &T_{s_1}\mapsto v, \quad T_{s_2} \mapsto v, \qquad 
%&\varepsilon \colon &T_{s_1}\mapsto -1, \quad T_{s_2} \mapsto -1, \\ 
%\varepsilon_1 \colon &T_{s_1}\mapsto v, \quad T_{s_2} \mapsto -v^{-1}, 
%\qquad &\varepsilon_2 \colon &T_{s_1}\mapsto -v^{-1}, \quad T_{s_2} 
%\mapsto v, \\ \rho\colon &T_{s_1} \mapsto \begin{bmatrix} -v^{-1} & 0 \\ 
%2v^{-1} &v \end{bmatrix}, \quad T_{s_2} \mapsto &\begin{bmatrix} v & -v\\ 0 
%& -v^{-1}\end{bmatrix}.& 
%\end{alignat*}

\begin{example} \label{expcanbas} 
Following Graham--Lehrer \cite{GrLe}, we define {\em cell representations} 
of $\cH$ as follows. Fix $\lambda \in \Lambda$ and let $W(\lambda)$ be a 
free $A$-module with basis $\{C_\fs \mid \fs\in M(\lambda)\}$. Then 
$W(\lambda)$ is an $\cH$-module with action given by 
\[ h.C_\fs=\sum_{\fs'\in M(\lambda)} r_h(\fs',\fs)\, C_{\fs'} \qquad 
\mbox{for $h \in \cH$ and $\fs \in M(\lambda)$}.\]
There is a symmetric bilinear form $g^\lambda \colon W(\lambda) \times 
W(\lambda) \rightarrow A$ defined by
\[ g^\lambda(C_\fs,C_\ft)=r_h(\fs,\fs) \qquad \mbox{where $\fs,\ft\in
M(\lambda)$ and $h=C_{\fs, \ft}^\lambda$}; \]
see \cite[Def.~2.3]{GrLe}.
We have $g^\lambda(h.C_\fs,C_\ft)=g^\lambda(C_\fs,h^*.C_\ft)$ for all 
$\fs,\ft \in M(\lambda)$ and $h \in \cH$; see \cite[Prop.~2.4]{GrLe}. 
Using the formula for $r_h(\fs',\fs)$ in the proof of 
Theorem~\ref{mainthm}, and a computation analogous to that in the proof 
of \cite[Prop.~20.6]{Lusztig03}, it is straightforward to check that 
\[v^{\ba_\lambda}\, g^\lambda(C_\fs,C_\ft)\equiv f_\lambda\,\beta_{\fs
\ft}^\lambda \quad \bmod vR[v].\]
In particular, we can see from this expression that the determinant
of the Gram matrix of $g^\lambda$ is non-zero. (This fact also follows 
from the general theory of cellular algebras, since $\cH_K$ is known to be 
split semisimple.)

Now let $K$ be the field of fractions of $A$; we write $\cH_K=
K\otimes_A \cH$ and $W_K(\lambda)=K\otimes_A W(\lambda)$.  Since
$g^\lambda$ induces a non-degenerate form on $W_K(\lambda)$ for 
each $\lambda \in \Lambda$, we have
\[ \Irr(\cH_K)=\{W_K(\lambda) \mid \lambda \in \Lambda\}; \qquad
\mbox{see \protect{\cite[Theorem~3.8]{GrLe}}}.\]
In fact, by the formula for $r_h(\fs',\fs)$ in the proof of 
Theorem~\ref{mainthm}, we see that the action of $\cH$ on $W(\lambda)$ is 
obtained by pulling back the action of $J$ on $E^\lambda_\spadesuit$ via 
Lusztig's homomorphism $\phi \colon \cH \rightarrow J_A$. This shows that 
$W_K(\lambda)\cong E_v^\lambda$; see \cite[20.2, 20.3]{Lusztig03}. 

Now let $\theta \colon A \rightarrow k$ be any ring homomorphism where $k$
is a field. By extension of scalars, we obtain a $k$-algebra $\cH_k=k 
\otimes_A \cH$ and cell representations $W_k(\lambda)=k \otimes_A 
W(\lambda)$ ($\lambda \in \Lambda$), which may no longer be
irreducible. Denoting by $g_k^\lambda$ the induced bilinear form on
$W_k(\lambda)$, we set $L^\lambda=W_k(\lambda)/\mbox{rad}(g_k^\lambda)$.
Then,  by \cite[Theorem~3.4]{GrLe}, each $L^\lambda$ is either $0$ or an 
absolutely irreducible representation of $\cH_k$, and we have 
\[ \Irr(\cH_k)=\{L^\mu \mid \mu \in \Lambda^\circ\} \qquad \mbox{where}
\qquad \Lambda^\circ=\{\lambda \in \Lambda \mid L^\lambda \neq 0\}.\]
Denote by $[W_k(\lambda):L^\mu]$ the multiplicity of $L^\mu$ as a
composition factor of $W_k(\lambda)$. Then
\begin{equation*}
\left\{\begin{array}{l} \quad [W_k(\mu):L^\mu]=1 \quad \mbox{for any
$\mu\in \Lambda^\circ$},\\[1mm] \quad [W_k(\lambda):L^\mu]=0 \quad 
\mbox{unless $\lambda \preceq \mu$ (i.e., $\lambda=\mu$ or $\ba_\mu<
\ba_\lambda$)}; \end{array}\right.  \tag{$\Delta$}
\end{equation*}
see \cite[Prop.~3.6]{GrLe}. Thus, the decomposition matrix 
\[D=\bigl([W_k(\lambda):L^\mu]\bigr)_{\lambda \in \Lambda,\mu \in 
\Lambda^\circ}\]
has a lower unitriangular shape, if the rows and columns are ordered 
according to increasing $\ba$-value.

The relations ($\Delta$) show that the subset $\Lambda^\circ \subseteq
\Lambda$ defines a ``canonical basic set'' in the sense of 
\cite[Def.~4.13]{mylaus},  which itself is an axiomatization of the results
obtained earlier by Rouquier and the author \cite{mykl}, \cite{GeRo2}. Note
that the definition of ``canonical basic sets'' only requires the 
invariants $\ba_\lambda$ (with respect to $L$) and the general set-up of 
Brauer's theory of {\em decomposition numbers} for associative algebras. 
Explicit descriptions of these ``basic sets'' are now known in all 
cases; see the survey \cite{mylaus} and the references there. In \cite{Jac2},
\cite{GeJa}, the existence of such basic sets for type $B_n$ and any $L$
(where $k$ has characteristic $0$) has been established {\em without} 
assuming Lusztig's conjectures on Hecke algebras with unequal parameters.
\end{example}

\begin{example} \label{expbn} Assume that $W$ is of type $B_n$ ($n\geq 2$) 
with diagram and weight function given by
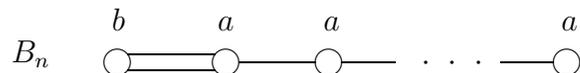
\begin{center}
\begin{picture}(250,30)
\put(  0, 10){$B_n$}
\put( 40, 10){\circle{10}}
\put( 44,  7){\line(1,0){33}}
\put( 44, 13){\line(1,0){33}}
\put( 81, 10){\circle{10}}
\put( 86, 10){\line(1,0){29}}
\put(120, 10){\circle{10}}
\put(125, 10){\line(1,0){20}}
\put(155,  7){$\cdot$}
\put(165,  7){$\cdot$}
\put(175,  7){$\cdot$}
\put(185, 10){\line(1,0){20}}
\put(210, 10){\circle{10}}
\put( 38, 22){$b$}
\put( 78, 22){$a$}
\put(118, 22){$a$}
\put(208, 22){$a$}
\end{picture}
\end{center}
where $a,b$ are non-negative integers. Then we have two cellular structures 
on $\cH$: the one given by Theorem~\ref{mainthm}, and the one given by 
Dipper--James--Murphy \cite{DJM3} (or, alternatively, by Graham--Lehrer 
\cite[\S 5]{GrLe}). In both cases, $\Lambda$ is the set of all pairs of 
partitions $(\lambda,\mu)$ such that $|\lambda|+|\mu|=n$. The partial order 
on $\Lambda$ that we use here is defined in terms of the $\ba$-function; it 
heavily depends on the weight function $L$, i.e., on the integers $a$ and $b$.
The partial order used by Dipper--James--Murphy (or Graham--Lehrer) is 
defined in purely combinatorial terms: it does not depend at all on the 
weight function $L$. Small examples show that the two partial orders indeed 
are different. Hence, in general, the cellular structure that we construct 
will be essentially different from the one constructed by 
Dipper--James--Murphy (or Graham--Lehrer). These differences can be seen 
explicitly in the descriptions of the sets $\Lambda^\circ \subseteq\Lambda$ 
parametrizing the irreducible representations of non-semisimple 
specialisations of $\cH$, as discussed in \cite{GeJa}. By 
\cite[Theorem~2.8]{GeJa}, the two parametrizations coincide if $b>(n-1)a>0$.
Thus, one may expect that the two cellular structures are equivalent under
this condition; see Iancu--Pallikaros \cite{IaPa} for a further discussion
of this question. 
\end{example}

\smallskip
\noindent {\bf Acknowledgements.} I wish to thank the referees for some useful
comments and for pointing out the reference \cite{kor2} to me.

%%%%%%%%%%%%%%%%%%%%%%%%%%%%%%%%%%%%%%%%%%%%%%%%%%%%%%%%%%%%%%%%%%%%%%%%%%%%%

\end{document}